\documentclass[a4paper, 11pt]{amsart}
\usepackage[latin1]{inputenc}
\usepackage[ T1]{fontenc}
\usepackage[english]{babel}
\usepackage{amssymb}
\usepackage{amsmath}
\usepackage{amsthm}
\usepackage{amscd}
\usepackage{amsfonts}
\usepackage{stmaryrd}
\usepackage{pb-diagram}
\usepackage{epic,eepic,epsfig}
\usepackage{a4wide}
\usepackage{nextpage}
\usepackage{fancyhdr}
\usepackage{enumerate}
\usepackage{hyperref}
\usepackage{float}
\usepackage{comment}
\usepackage{todo}
\usepackage{faktor}
\usepackage{listings}
\usepackage{diagbox}
\usepackage{dsfont}

\usepackage[
backend=biber,
style=alphabetic,
]{biblatex}
\newtheorem{teo}{Theorem}

\newtheorem{lema}[teo]{Lemma}
\newtheorem{prop}[teo]{Proposition}

\newtheorem{rem}[teo]{Remark}

\renewcommand{\Im}{\operatorname{Im}}

\newcommand{\Z}{\mathbb{Z}}
\newcommand{\Q}{\mathbb{Q}}
\newcommand{\Qbar}{\overline{\Q}}

\newcommand{\N}{\mathbb{Z}_{> 0}}
\newcommand{\C}{\mathbb{C}}

\newcommand{\eps}{\varepsilon}

\newcommand{\D}{\mathds{D}}

\newcommand{\al}{\alpha}
\newcommand{\be}{\beta}

\renewcommand{\H}{\mathcal{H}}

\newcommand{\SL}{\operatorname{SL}}

\newcommand{\trdeg}{\operatorname{trdeg}}

\begin{document}

\title{The St\'ephanois theorem with only prime isogenies}
\author{Desir\'ee Gij\'on G\'omez}

\address{Department of Mathematical Sciences, University of Copenhagen,
Universitetsparken 5, 
2100 Copenhagen O, Denmark.}
\email{dgg@math.ku.dk}

\curraddr{INRIA Nancy-Grand Est, Campus Scientifique, BP 239
54506 Vandoeuvre-l\'es-Nancy, France}
\email{desiree.gijon-gomez@inria.fr}

\thanks{We thank Fabien Pazuki for guidance in this project and for encouraging us to submit it. This material is part of the author's PhD thesis, and was written while being a PhD student at the University of Copenhagen, with minor corrections carried out during her current position as a postdoc in INRIA Nancy-Grand Est.}

\thispagestyle{empty}

\begin{abstract}
We present a strengthening of the proof of the St\'ephanois theorem \cite{Stephanois}. We follow the modular version by \cite[Th\'eor\`eme 1]{WaldschModVerStephanois},\footnote{which is based in a suggestion by Daniel Bertrand \cite{BertThetaFuncTransc}} but it also applies to the original proof. The improvement is not in the result or the conditions, but in the need of weaker tools on the proof itself. More precisely, we only employ modular polynomials of \emph{prime} degree, instead of polynomials of arbitrary level. Furthermore, one can restrict to primes in fixed arithmetic sequence. 

On the proof itself, the only crucial difference appears in \cite[Cinqui\`eme pas, page 5]{WaldschModVerStephanois} proof, and on the final contradiction in \cite[Septi\`eme pas, page 7]{WaldschModVerStephanois}, but for readability, we present a complete proof with this modification. 

This is part of a larger project to generalize the St\'ephanois theorem to the Igusa invariants of curves of genus two, as the Siegel modular polynomials in the literature are usually only considered for prime levels. 

The material is part of Chapter 7 of the author's PhD thesis.
\end{abstract}

 \maketitle
\begin{center}
---------
\end{center}
\tableofcontents

\section{Introduction}
Consider the elliptic $j$-invariant $j: \H \to \C$, a $\SL_2(\Z)$-invariant modular function, which classifies elliptic curves over $\Qbar$. It admits a Fourier expansion with integer coefficients:
\[
J(q) = \frac{1}{q} + 744 + 196884q + \cdots \, q = e^{2\pi i \tau},
\]
which defines a meromorphic function on $\D = \{ q \in \C:$ $|q| < 1\}$, the only pole being at $q=0$ and simple. This short note is about the following result. 

\begin{teo}[The St\'ephanois\footnote{It was previously called the Mahler-Manin conjecture. The French name of this result is \emph{le t\'eor\`eme st\'ephanois}. We claimed in \cite{OnCMException} that we could not find an English translation to this result, but we found one in \cite{Pellarin2011arxiv} and \cite{AdamMahlerMethod}, with the English translation being \emph{Stephanese}. The explanation of this name is that it was proven by a team of four people in Saint-\'Etienne, and in French \emph{st\'ephanois} is the gentilic for that region in France.} theorem] Let $q \in \D \setminus \{0\}$, then $\trdeg \Q(q, J(q)) \geq 1$.
\end{teo}

It can be argued that the proof of \cite{Stephanois} is an adaptation of what is called the \emph{Mahler method}. For more references on that matter and the Mahler method in the context of the modular invariant, see \cite{Gramain}, \cite{PellarinMahlermethod} and \cite{AdamMahlerMethod}. 

We are going to present here the "modular" proof of the St\'ephanois theorem (\cite{Stephanois}) in \cite{WaldschModVerStephanois}. The main difference in both approaches lies on how to cancel out the pole at infinity of the $j$-function. The original one was simply to consider $\tilde{J}(q):= q J(q)$, and the modular approach uses the modular discriminant form $\Delta$, as it has a simple zero at infinity, and $\Delta(q) J(q)$. That implies the following differences.
\begin{itemize}
\item The upper bounds of the auxiliary function $F(q)$ rely on upper bound on the (integral) Fourier coefficients of either $\tilde{J}^l(q)$ or $\Delta^k J(q)^l$ for $k \geq l$. In the first case, they use Mahler's bound on the coefficients of powers of the $j$-invariant from \cite[Equation (3), page 201]{MahlerOnTheCoeffOf}. In the modular proof, one observes that, \emph{independently of} $l$, $\Delta^k J(q)^l$ is a \emph{cusp} form of weight $12 k$, and one can invoke the classical Hecke's bound on the Fourier coefficients of cusp forms. 
\item There is a step requiring lower bounds on the auxiliary function $F$ on $q^S$ for a fixed $S \in \N$, where $q$ is an (eventually impossible) point in $\D$ such that $(q, J(q)) \in \Qbar^2$. On the original proof, by construction $F(q) = A(q^S, J(q^S))$ is a polynomial with integer coefficients evaluated in the pair of algebraic numbers $(q^S, J(q^S))$. 

Then the lower bound follows from Liouville's inequality on the algebraic number $F(q^S, J(q^S))$. In the modular proof, the auxiliary function is now $F(q) = \Delta^{2N} A(q, J(q))$ (where $\deg_X A< N$ and $\deg_Y A < N$). $\Delta^{2N}(q^S)$ is not algebraic anymore, but one combines the previous lower bound on $\Delta^{-2N}(q^S) F(q^S) = A(q^S, J(q^S)) \in \Qbar$, and a lower bound on $\Delta^{2N}(q^S)$ that is negligible in comparison with the first one.
\end{itemize}

The main strategy (in both the original and modular proof) goes back to the Mahler method: one starts with a pair $(q, J(q)) \in \Qbar^2$. For any $n \in \N$, then $(q^n,  J(q^n)) \in \Qbar^2$, because $J(q^n)$ is algebraic over $\Q(J(q))$ via the modular polynomial $\Phi_n$ of level $n$. For the proof, one requires further information on the modular polynomial, namely upper bounds on its degree and height. 

Our original contribution is that it is possible to restrict oneself to only work with \emph{prime} powers $\{q^p\}_{p \text{ prime}}$,\footnote{an unavoidable unfortunate notation, $q$ is the complex variable $q = e^{i2\pi \tau}.$} hence only working with isogenies of prime level, and only use bounds for the modular polynomials of prime level. Following \cite[Th\'eor\`eme 1]{WaldschModVerStephanois}, the only modifications come at Cinqui\`eme pas, on the determination on the fixed power $P$, and Septi\`eme pas, on the final contradiction in the transcendence proof. This restriction gives a worse bound on $P$, but it is still enough to conclude in the contradiction. 

Morally, one should expect this modification to work, because every isogeny between elliptic curves factors as composition of prime isogenies. In a way, all the information encoded in the modular polynomials in every level is already encoded in the modular polynomials of prime level. 

Furthermore, it also follows from the arithmetic bounds that one can restrict even more to work with primes in fixed arithmetic sequences.

This is part of a work to generalize this result to the Igusa invariants for genus two curves in the author's PhD thesis. The modular proof is more suitable for this goal: the corresponding invariants are also defined as quotients of Siegel modular forms, with the denominator being (powers of) a cusp form, and its zero locus is the only subset where the invariants are not defined. Hence, as it happens with $\Delta$, it can be used to cancel out the denominators. More importantly, Hecke's bound on Fourier coefficients on cusp forms generalizes to Siegel cusp forms. Our contributions with the restriction to prime level isogenies is because the modular polynomials for genus two curves (Siegel and Hilbert modular polynomials) are usually only considered for prime levels in the literature.   

\subsection{List of parameters}

As it happens with a transcendence proof, we work with integral parameters, dependent on each other, such that we have an impossible inequality in the end. Alongside, there will be an array of constants (some absolute, others dependent only on $q$) appearing throughout the proof, for which we reserve the notation $C_i$ for $i=0,1,\dots$. We start the proof with a fixed $q \in \D$ such that $(q,J(q)) \in \Qbar$. For the complex variable in $\D$, we set $z$.
\begin{description}
\item $N \in \N$ will appear first as $\deg_X A< N$ and $\deg_Y A < N$, for $A \in \Z[X,Y]$ the auxiliary polynomial. All other functions and parameters depend ultimately and explicitly on $N$. Its sole purpose is being large enough for the contradiction to happen in the end. It depends on $|q|$ only in Equation \eqref{EqDependNonq} via a lower bound for $N.$
\item $L \in \N$ will appear in the application of the Siegel's lemma for $F(z) = \Delta^{2N} A(z, J(z))$. The polynomial $A$ is constructed so that $F$ vanishes at $z=0$ with order at least $L$. In particular, one should have $N^2 \geq 2L$ for the bounds of the Siegel's lemma. It will be set $L := [N^2/2]$, for $[ \cdot ]$ the integer part, therefore $L$ and $N$ are related via
\[
2L \leq N^2 < 2(L+1).
\]
\item $M \in \N$ will be defined as the actual order of vanishing of the auxiliary function $F$ at $0$, hence 
\[
M \geq L.
\]
\item $P \in \N$ will be the prime number such that we eventually consider $q^P$. It is the only parameter depending more intrinsically on $q$. In the standard Mahler's method, this parameter tends to infinity, but in this proof it cannot be arbitrarily large. Hence it will be defined as the minimum $P$ possible, and one needs it to not be too large in terms of the other parameters for the contradiction to happen. More precisely, in \eqref{EqBoundOnS} we will prove
\[
\frac{P^2}{\log P} = O_{|q|}\left(N\log M \right).
\]
\end{description}

\section{Preliminary lemmas}
\label{sec:PreliminaryModularStep}
Here we state the preliminary necessary results.

\begin{prop} The meromorphic function $J: \D \to \C$ induced by the Fourier expansion of the $j$-invariant is transcendental over $\Q(z).$ \label{PropTranscendenceOfJ}
\end{prop}
\begin{proof} See \cite[Lemme 4]{Stephanois}, or \cite{MahlerAlgDiffEquation}, that we include for completeness.  Equivalently, we are going to prove that, in the $\tau$ variable, the functions $j: \H \to \C$ and $q(\tau) = e^{i2\pi\tau}$ are algebraically independent over $\Q$, using that $j$ is $\SL_2(\Z)$-invariant.

Assume then there exists a polynomial $P \in \Q[X,Y]$ such that, as function in $\H$, $P(q,j) \equiv 0$, or
\begin{equation}
p_n(j) q^n + \ldots + p_0(j) = 0, \label{EqFunctionalEquationqj}
\end{equation}
for $p_i \in \Q[X]$, $p_n \not\equiv 0$. For $\gamma = \begin{pmatrix}
a& b \\ c & d
\end{pmatrix} \in \SL_2(\Z)$ set $\gamma \tau = \frac{a\tau + b}{c\tau + d}$ the corresponding linear fractional transformation. 

Consider $\tilde{\tau}$ such that $p_n(j(\tilde{\tau})) \not= 0$ and the polynomial
\[
P_{\tilde{\tau}}(X) =  X^n + \frac{p_{n-1}(j(\tilde{\tau}))}{p_n(j(\tilde{\tau}))}X^{n-1} + \ldots + \frac{p_0(\tilde{\tau})}{p_{n}(j(\tilde{\tau}))}.
\] 
Then $P(q, j) \equiv 0$ implies $P(q\circ \gamma, j\circ \gamma) \equiv 0$ for any $\gamma \in \SL_2(\Z)$, and as $j \circ \gamma =j$, then $P_{\tilde{\tau}}$ vanishes on $\{ q(\gamma \tilde{\tau}),$ for all $\gamma \in \SL_2(\Z)\}$. 

It is then enough to check that $\{ q(\gamma \tilde{\tau})\}_{\gamma \in \SL_2(\Z)}$ has infinitely many distinct elements. The argument from \cite{Stephanois} takes $\gamma_c = \begin{pmatrix}
1 & 0 \\ c & 1
\end{pmatrix}$, for $c \in \N$ and $\tilde{\tau} = it$ with $t \geq 1$ a transcendental number. Then if $c\not= c'$, $q(\gamma_c(it))\not= q(\gamma_{c'}(it)$, for that they can only be equal if $\gamma_c(it) - \gamma_{c'}(it) = k$ for some $k\in \Z$, which is equivalent to $(c-c' + kcc')t^2 - ik(c+c')t - k = 0$, which as $t$ is transcendental implies $k = c-c' = 0$. 

Therefore $p_l(j(it)) = 0$ for all $t>1$ transcendental, and as $j$ attains all its real values on the vertical line, necessarily $p_l \equiv 0$ for $l=1, \ldots n-1$, and hence $P \equiv 0.$

An alternative argument comes from \cite{MahlerAlgDiffEquation}. For a rational number $\theta \in \Q$, there exists a sequence $\gamma_n = \begin{pmatrix}
a_n & b_n \\ c_n & d_n
\end{pmatrix} \in \SL_2(\Z)$ with 
\begin{description}
\item $\frac{a_n}{c_n} \to \theta$,
\item $\frac{c_n}{d_n} \to 0$, $|c_n d_n| \to \infty$, with $d_n \geq 0$. Note then that $|c_n(c_n \tau + d_n)| \to \infty$ for any $\tau \in \H.$
\end{description}
As $\gamma \tau  = \frac{a\tau + b}{c\tau + d}= \frac{a}{c} - \frac{1}{c(c(\tau + d)}$, then 
\[
e^{2\pi i \gamma_n \tau} = e^{2\pi i \frac{a_n}{b_n}} e^{-2\pi i \frac{1}{c_n(c_n \tau + d_n)}} \to e^{2\pi i \theta}, \text{ when } n \to \infty.
\]
Taking limits in \eqref{EqFunctionalEquationqj}, it follows that the set of roots of $P_{\tilde{\tau}}$ contains infinitely many roots of unity, and we also arrive at the same conclusion.
\end{proof}

\begin{rem} In \cite{MahlerAlgDiffEquation}, it is furthermore proven that the four functions \\$q(\tau) = e^{i2\pi \tau}, j(\tau), j'(\tau), j''(\tau)$ are algebraically independent over $\C(\tau)$ (and this is best possible as $j$ solves a differential equation of order three). This result was extended to any dimension in \cite{BerZuDTransDeg} and \cite{BerZud}.
\end{rem}

As we said before, we rely on results for the Fourier coefficients of the Fourier expansion of a holomorphic cusp form. 

\begin{lema}[Hecke bound on cusp forms] There exists a constant $C_1 >0$ such that the following holds. For any $N \in \N$ and any $1 \leq l \leq N$, the functions $\Delta^{2N} J^l$ are modular cusp forms of weight $24 N$, and if we consider their Fourier expansions at infinity:
\[
\Delta^{2N} J^l(q) := \sum_{k \geq 1} c_{N,l}(k) q^k,
\]
then $c_{N,l}(k) \in \Z$ for all $N,k,l$ and
\[
|c_{N,l}(k)| \leq C_1^N k^{12N}.
\] \label{LemmaHeckeBound}
\end{lema}
\begin{proof}
First $\Delta^{2N} J^l$ is holomorphic and vanishes at infinity, for that $\Delta$ has a simple zero at infinity that cancels out with the pole of $J$ for our choices of parameters $l \leq N$. Integrality of the Fourier coefficients follows from $\Delta$ and $J$ having integral Fourier coefficients. It is a modular form of weight $24N,$ and by the classical Hecke estimate for cusp forms, $|c_{N,l}(k)| = O(k^{12N})$.

More precisely, we have a convenient choice of explicit constant, from the proof of the Hecke estimate (see \cite[Chapter VII, Section 4, Theorem 5]{SerreCourseArithm}). For a cusp form $f = \sum_{k \geq 1} f(k) q^k$ of weight $\omega$ then the function $\phi_f(\tau) := \Im(\tau)^{\omega/2}f(\tau)$ is bounded in $\H$, and setting $	C(f):= e^{2\pi} \sup_{\H} |\phi_f|$, then $|f(k)| \leq C(f)k^{\omega/2}$ for all $k\geq 1$. Observe that then $C(gh) \leq C(g)C(h)$ for any cusp forms $g$ and $h$, hence we can deduce explicit constants for the forms $C(\Delta^{2N} J^l)$. Set $C_1 := \max \{ 1, C(\Delta^2), C(\Delta^2 J) \}$, then for any $N$ and $1 \leq l \leq N$,
\[
C(\Delta^2N J^l) = C\left(\Delta^{2(N-l)} (\Delta^2 J)^l \right) \leq C_1^{N-l + l} = C_1^N,
\]
and the result follows.
\end{proof}

The auxiliary function to consider in the proof is a polynomial \emph{with integer coefficients} that solves a suitable linear system over $\Z$. There is a standard result to use in this situation, known as Siegel's lemma. Roughly, for an homogeneous system of $X$ variables and $Y$ equations with coefficients in $\Z$, and if we allow $X >Y$, then the kernel is positive dimensional as a vector space, and we can find \emph{integral} solutions in that kernel. Siegel's lemma gives one such solution that it is small in terms of the size of the system.
\begin{prop}[Siegel's lemma] Consider the linear system with coefficients in $\Z$ in $X$ variables and $Y$ equations with $X > Y$:
\[
\sum_{j=1}^{X} m_{i,j}x_{i,j} = 0, \quad 1 \leq j \leq Y,
\]
Set $B = \max |m_{i,j}|$, then there exists a non zero solution in $\Z^X$ such that 
\[
|x_{i,j}| \leq \left(XB\right)^{\frac{Y}{X-Y}}, \label{PropSiegelsLemma}
\]
\end{prop}
A standard extra assumption for applying Siegel's lemma is to have $X \geq 2Y$, so that $Y/(X-Y) \leq 1$.  

\subsection{Height lemmas}
We also state Liouville's inequality\footnote{There appears to be various results in Diophantine approximation with this name.} for algebraic numbers, and another results of height nature. We use the following notations for $\al \in \Qbar$:
\begin{description}
\item $\deg(\al)$ for its degree, i.e.\ $\deg(\al) = [\Q(\al): \Q]$;
\item $h(\al)$ for its logarithmic (Weil) height;
\item $M(\al)$, for $\al \in \Qbar$ for the Mahler measure of the minimal polynomial of $\al$ in $\Z[X];$
\item $m(\al) = \log M(\al)$ for its logarithmic Mahler measure; remark then that $m(\al) = \deg(\al) h(\al)$;
\item for any polynomial $P \in \Z[X_1, \ldots X_r]$ we denote its \emph{height} and its \emph{length} as $H(P) := \max |c|$ and $L(P):= \sum |c|$ with $c$ ranging over the coefficients of $P$. Likewise, we have lowercase notation for its logarithms: $h(P) : = \log h(P)$, and $l(P) = \log L(P).$
\end{description}

\begin{lema}[Liouville's inequality] Let $0 \not= \al \in \Qbar$, then it holds that
\[
\log |\al| \geq -\deg(\al)h(\al).
\] \label{LemmaLiouvilleIneq}
\end{lema}
\begin{proof}
    See \cite[Theorem 11.5.21]{BombieriHeights}.
\end{proof}

The following lemma allows us, for $\al, \be \in \Qbar$ linked via $P(\al, \be) = 0$ for $P \in \Z[X,Y]$, to control arithmetic information of $\be$ in terms of arithmetic information of $\al$ and $P$. It will be used for the modular polynomials.

\begin{lema} Consider $\al, \be \in \Qbar$ and a polynomial $P \in \Z[X,Y]$ such that $P(\al, \be) = 0$. Assume that $P(\al, Y)$ is not a constant polynomial. Then
\[
m(\be) \leq \deg(\al) \left( \log L(P) + \deg_x(P) h(\al) \right).
\]\label{LemaMahlerMeasOnARoot}
\end{lema}
\begin{proof}
See \cite[Lemme 5]{Stephanois}.
\end{proof}

Finally, we are going to require bounds of the height of the evaluation of an integral polynomial in algebraic values. Remark that this is a generalization of the standard properties of the height for $\al, \be \in \Qbar$:
\[
h(\al\be)  \leq h(\al) + h(\be), \text{ and }h(\al \pm \be)  \leq \log 2 + h(\al) + h(\be).
\]

\begin{lema} Consider $P \in \Z[X_1, \ldots X_n]$ a non zero polynomial, and  $\al_1, \ldots \al_n \in \Qbar$. Then the following holds
\[
h(P(\al_1, \ldots , \al_n)) \leq \log L(P) + \sum_{i=1}^n \deg_{X_i}(P) h(\al_i).
\]  \label{LemaHeightEvalPol}
\end{lema}
\begin{proof}
    See \cite[Lemma 3.7]{WaldscDiophAproxLinAlg}.
\end{proof}
This result, together with the properties of the the modular polynomials from the following subsection, was the original approach in \cite{Stephanois} for comparing $h(J(q^p))$ with $h(J(q))$. Alternatively, there is a more geometric approach in \cite[Lemma 2]{BertThetaFuncTransc}, by passing through the Faltings height.

\begin{prop} There exists $C_2>0$ such that for any $q \in \D$ such that $J(q) \in \Qbar$, and any $n \in \Z,$
\[
h(J(q^n)) \leq 2h(J(q)) + 6\log(1 + n) + C_2.
\] \label{PropHeightIsogeJInv}
\end{prop}
\begin{proof}
See \cite[Lemma 2]{BertThetaFuncTransc}.
\end{proof}

\subsection{Modular polynomials}
Another main character are the modular polynomials $\Phi_N(X,Y) \in \Z[X,Y]$ for $N \geq 1$. They are characterized for being the (irreducible and monic) polynomial in two variables such that $\Phi_N(j_1, j_2)=0$ if and only if $j_1$ and $j_2$ are the $j$-invariants of elliptic curves which are isogenous via a cyclic isogeny of degree $N$. We have the following results on its degree and height from \cite{PublishExplicitBounds}.
\begin{lema} Consider the modular polynomial $\Phi_N(X,Y) \in \Z[X,Y]$ for a \emph{general} level $N$.
\begin{itemize}
\item Its degree $\deg(\Phi_N):= \deg_X(\Phi_N) = \deg_Y(\Phi_N) = \psi(N)$ where $\psi$ is the Dedekind psi function given by 
\[
\psi(N) = N \prod_{p \lvert N \text{ prime }} \left( 1 + \frac{1}{p} \right).
\]
\item\emph{(Cohen's bound)} It follows that
\[
h(\phi_N) = 6 \psi(N)( \log N - 2 \kappa_N + O(1) ),
\]
where $\kappa_N = \sum_{p \lvert N \text{ prime }} \frac{\log p }{p}$.
\item \emph{(Explicit bound)} If one replaces $\kappa_N$ with $\lambda_N = \sum_{p^n \lvert \lvert N} \frac{p^n-1}{p^{n-1}(p^2 -1)}\log p$, then one can take the explicit constants:
\[
6 \psi(N) \left( \log N - 2 \lambda_N - 0.0351 \right) \leq h(\Phi_N) \leq 6 \psi(N) \left( \log N - 2 \lambda_N + 9.5387 \right).
\]
\end{itemize}
\end{lema}

We state the best results known so far, but for the sake of the proof, Mahler's bound $h(\phi_N) = O(N^{3/2})$ from \cite{MahlerOnTheCoeffOf} would have been enough. We specialize for prime level, for future use. The explicit bound in this case is from \cite{BrSu}.

\begin{lema} Consider $\Phi_p(X,Y) \in \Z[X,Y]$ for $p$ a prime number. 
\begin{itemize}
\item Its degree $\deg (\phi_p) = p+1$.
\item We have that
\[
h(\phi_p) = 6p\log p + O(p).
\]
\item (Explicit bound) We have that
\[
h(\phi_p) \leq 6p\log p + 16p + 14\sqrt{p}\log p.
\]
\end{itemize} \label{LemaModPolPrime}
\end{lema}
Therefore, we have the bound on $L(\phi_p) \leq (\psi(p)+1)^2 e^{h(\phi_p)}$, so 
\[
\log L (\phi_p) \leq 2\log(p+2) + 6p\log p + 16p + 14\sqrt{p}\log p = O(p \log p)
\]

Finally, we state the following result about the degree of the specialization of the modular polynomials, from \cite[lemma 2 ii)]{BertThetaFuncTransc}.

\begin{prop} Let $q \in \D$ such that $J(q) \in \Qbar$, and let $E$ the elliptic curve with $J(E) = J(q)$. Then
\begin{itemize}
\item If $E$ is CM, consider $\tau \in \H$ a period for $E$ and equation $a\tau^2 + b\tau + c =0$. Then if $p$ does not divide $a$:
\[
\frac{p-1}{3} \leq [\Q(J(q),J(q^p)): \Q(J(q))].
\]
\item If $E$ is \emph{not} CM, consider $p_E$ the prime number given by Serre's theorem in \cite[Th\'eor\`eme 2]{SerreOpenImageTheorem}. Then if $p > p_E$,
\[
p+1 = [\Q(J(q),J(q^p)): \Q(J(q))].
\]
Equality to $p+1$ also holds in the first case for $p$ being inert in the CM field of $E$.
\end{itemize} \label{PropBertrand}
\end{prop}

\subsection{Bounds on sums of primes}
This last subsection is for stating some results of classical analytic number theory, as consequences of the prime number theorem. 

\begin{prop} Let $\pi(x) := \#\{p \text{ prime, } p \leq x\}$ the prime counting function.
\begin{itemize}
\item For $x \geq 2$,
\[
\pi(x)  = \frac{x}{\log x} + \frac{x}{(\log x)^2} + O\left(\frac{x}{(\log x)^3}\right).
\]
But we may use the following weaker version.\footnote{sometimes called Tchebychev's theorem.} There exists an absolute constant $C_{14}>0$ such that  
\[
\frac{x}{\log x} \leq \pi(x) \leq C_{14}\frac{x}{\log x}. 
\]
\item For $x$ large enough, 
\[
\frac{x^2}{2\log x} \leq \sum_{p \text{ prime, } p < x} p \leq \frac{x^2}{\log x}.
\]
\end{itemize}

Let $\delta \in \N$ and $a \in \N$ such that $\gcd(a,\delta) = 1$, and consider 
\[
\pi(x,a,\delta) : =\# \{ p \text{ prime, } p\equiv a \text{ (mod }\delta) \}.
\]
\begin{itemize}
\item For $x$ large enough only in terms of $\delta$, 
\[
\pi(x)  = \frac{x}{\phi(\delta)\log x} + \frac{x}{\phi(\delta)(\log x)^2} + O_{\delta}\left(\frac{x}{(\log x)^3}\right),
\]
for $\phi$ Euler's totient function.
\item
For $x$ large enough only in terms of $\delta$.
\[
\frac{x^2}{\phi(\delta)2\log x} \leq \sum_{p \text{ prime, } p\equiv a(\delta), p < x} p \leq \frac{x^2}{\phi(\delta)\log x}.
\]
\end{itemize} \label{PropBoundsSumsPrimes}
\end{prop}
\begin{proof}
The first statement is the quantitative version of the prime number theorem in \cite[Chapter 6, equation under (6.15)]{MontVauMNT}. For the second statement,\footnote{we saw this \href{https://mathoverflow.net/questions/63412/upper-bounds-for-the-sum-of-primes-up-to-n}{here}.} set $\Sigma(x) = \sum_{p \text{ prime, } p < x}p$. We require the prime number theorem until the second order term solely to have cleaner constants in our bounds. By Abel's summation formula and the prime number theorem above,
\[
\Sigma(x) = \pi(x)x - \int_{t=2}^x \pi(t) dt = \frac{x^2}{\log x} + \frac{x^2}{(\log x)^2} + O\left(\frac{x}{(\log x)^3}\right) - \int_{t=2}^x \pi(t) dt.
\]
For any $m \geq 1$, by integration by parts, \small
\[
\int_{2}^x \frac{t}{(\log t)^m} dt = \int_{2}^x \frac{1}{(\log t)^m} t dt = \left.\frac{t^2}{2(\log t)^m}\right\rvert_{2}^x + \int_{2}^x \frac{t}{m(\log t)^{m+1}} dt = \frac{x^2}{2(\log x)^m} + O \left(\frac{x^2}{m(\log x)^{m+1}} \right),
\]
\normalsize hence \small
\begin{align*}
\int_{t=2}^{x} \left( \frac{t}{\log t} + \frac{t}{(\log t)^2} \right) dt = \frac{x^2}{2 \log x} + \frac{3}{2}\int_{t=2}^x \frac{t}{(\log t)^2} dt + O(1) = \frac{x^2}{2 \log x}  + \frac{3}{4}\frac{x^2}{(\log x)^2} + O\left( \frac{x^2}{(\log x)^3} \right).
\end{align*}
\normalsize so 
\[
\int_{t=2}^{x} \pi(t) dt = \frac{x^2}{2 \log x} + \frac{3}{4}\frac{x^2}{(\log x)^2} + O\left( \frac{x^2}{\log x^3} \right),
\]
and it follows
\[
\Sigma (x) = \frac{x^2}{2 \log x} + \frac{x^2}{4(\log x)^2} + O\left(\frac{x^2}{(\log x)^3} \right).
\]
And one can check computationally with SageMath \cite{Sage} that already for $x \geq 11$,
\[
\frac{x^2}{2\log x} \leq \sum_{p \text{ prime, } p < x} p \leq \frac{x^2}{\log x}.
\]
For the statements on primes in arithmetic progression, the relevant statement is \cite[Corollary 11.21]{MontVauMNT}. Also, for more explicit versions see \cite[Theorem 1.1]{BoundExplicitPrimes}
\end{proof}

\section{Construction of the auxiliary function}
Our auxiliary function is a polynomial $A \in \Z[X,Y]$ with $\deg_X A < N$ and $\deg_Y A < N$ such that the auxiliary function $F(z) = \Delta(z)^{2N}A(z,J(z))$ is holomorphic and with high vanishing order at $z=0$, say at least $L$. If we write $A(x,y) = \sum_{0 \leq i,l < N} a_{i,l} x^i y^l$, we can consider the Fourier expansion of $\Delta(z)^{2N}A(z,J(z))$ and set the first $L$ of them to be zero. More precisely, for $\Delta^{2N}J^l=\sum_{k \geq 1} c_{N,l}(k) z^k$ as in Lemma \ref{LemmaHeckeBound}, we write
\begin{align}
\begin{split}
\Delta(z)^{2N}A(z,J(z)) &= \sum_{0 \leq i,l < N} a_{i,l} z^i \left(\Delta(z)^{2N} J(z)^{l}\right) = \sum_{0 \leq i,l < N} a_{i,l} z^i \sum_{k \geq 1} c_{N,l}(k) z^k = \\
&= \sum_{0 \leq i,l <N} \sum_{k \geq 1} a_{i,l} c_{N,l}(k) z^{i + k}
= \sum_{\nu \geq 0} \left(\sum_{\substack{0 \leq i \leq \min{\nu, N-1} \\ 0 \leq L < N}} a_{i,l} c_{N,l}(\nu -i) \right)z^\nu.
\end{split} \label{EqTaylorExpansionF}
\end{align}
We consider the system in $N^2$ variables and $L$ equations,
\[
\sum_{0=i}^{\min{\nu, N-1}} \sum_{l=0}^{N-1}  c_{N,l}(\nu -i) a_{i,l}, \text{ for } 0 \leq \nu < L,
\]
and we apply Siegel's lemma Proposition \ref{PropSiegelsLemma}, under $N^2 \geq 2L$. The upper bound on the coefficients on the system comes from Lemma \ref{LemmaHeckeBound}, that gives us $|c_{N,l}(\nu - i)| \leq C_1^N (\nu - i)^{12N}$, so $B \leq C_1^N L^{12N}$, and hence
\begin{equation}
\begin{split}
|a_{i,l}| & \leq N^2 C_1^N L^{12N}, \\
L(A) &:= \sum_{0 \leq i,l <N} |a_{i,l}| \leq N^4C_1^N L^{12N} \leq_{\exists C_3} C_3^N  L^{12N}.
\end{split} \label{EqSielegLemmaAuxiliaryPolynomial}
\end{equation}

Finally, $F$ is not the constantly zero function by Proposition \ref{PropTranscendenceOfJ}.

\section{Upper bound on $|F(z)|$}
Set $M = ord_{z=0}F(z)$, for $F$ our auxiliary function. By construction, $M \geq L$. The following is an upper bound on $F$ on closed disks for $0<r<1$, $\D_r \subset \D_1$ of the form $|F(z)| = O_{M,N}(z^M)$. As $z^{-M}F(z)$ is holomorphic on $\D_1$, by definition of $M$, we could apply the maximum principle to $G(z):= z^{-M}F(z)$ on closed disks\footnote{The boundary of $\D_1$ is the natural boundary of modular forms, so these arguments need to be carried out on $\D_r$.} $\D_r$, but then we would need the the dependence of $\max _{\D_r} G$ in terms of $N,M,r$. Instead, one uses that the Taylor expansion of $G(z) = \sum d_{k} z^{k}$ is of polynomial growth, meaning the terms $d_{k} = O(k^{rN})$ for $r \in \N$, so they are comparable to derivatives and powers of $\frac{1}{1-z} = \sum_{n \geq 0} z^n$. 

With more details, by \eqref{EqTaylorExpansionF},
\[
G(z) := z^{-M} F(z) = \sum_{k \geq 0} d_{k} z^{k},
\]
for 
\[
d_{k} := \sum_{0 \leq i,l < N} a_{i,l} c_{N,l}(M + k -i), 
\]
where the sum has $N^2$ terms,\footnote{Remark that the restriction $i < k$ from \eqref{EqTaylorExpansionF} does not apply anymore.} $|c_{N,l}(M + k -i )| \leq C_1^N(M + k -i)^{12N}$ by Lemma \ref{LemmaHeckeBound} and $|a_{i,l}| \leq N^2C_1^NL^{12N}$ by \eqref{EqSielegLemmaAuxiliaryPolynomial}, therefore exists $C_4 >0$ with
\[
|d_k| \leq N^4C_1^NL^{12N}C_1^N (M +k)^{12N} \leq C_4^N L^{12N} (M+k)^{12N}.
\]
The following is a pretty general argument for power expansion with coefficient of this type of growth: bound a function $\sum_{t=1}^{\infty} (M + t)^{12N}$ via comparing it with the expansion of $1/(1-x)$, powers of it, or derivatives.

\begin{lema}[Explicit Schwarz lemma for integral series of polynomial growth] Assume \\$h(z) = \sum_{t \geq 0} h_n z^n  \in \Z[[z]]$ is a holomorphic function on $\D$ such that $|h_n| \leq (M + n)^{K}$ for fixed $K \in \N$. Then
\[
|h(z)| \leq (M+1)^{K} K! \frac{1}{(1 - |z|)^{K+1}}. \label{LemaExplicitSchwarz}
\]
\end{lema}
\begin{proof}(From \cite[proof of Lemma 2.1]{Netserenko})
Note that the Taylor expansion for any $k \geq 1$ and any $0 < x <1,$
\[
\frac{k!}{(1-x)^{k +1}} = \sum_{n=0}^{\infty} (n+1)\cdots (n+k) x^k,
\]
hence 
\[
1 + \sum_{n\geq 1}n^k x^{k} \leq \sum_{n \geq 0} (n+1)\cdots (n+k) x^k = \frac{k!}{(1-x)^{k +1}},
\]
so we have the bounds (remark $(M + n)^{K} \leq (M+1)^{K} n^{K}$):
\[
\sum_{n\geq 0} (M + n)^{K}|z|^n \leq (M+1)^{K} \left( 1 + \sum_{n\geq 1} n^{K}|z|^n \right) \leq (M+1)^{K} \frac{K!}{(1 - |z|)^{K+1}}.
\]
\end{proof}

By Lemma \ref{LemaExplicitSchwarz} above, it follows that 
\[
|F(z)| \leq |z|^M C_5^N L^{12N}M^{12N}N^{12N}\frac{1}{(1 - |z|)^{12N + 1}}.
\]
We establish now the only dependence of $N$ on $q$, setting $N$ large enough so that for $r = \frac{1 + |q|}{2}$,\footnote{This is simply so one argues on a disk containing $q$ in its interior, and this is an easy enough choice of radius.}
\begin{equation}
\left(\frac{1}{(1 - r)}\right)^{12N + 1} \leq \left( \frac{N^2}{2} \right)^{N}, \label{EqDependNonq}
\end{equation}
that in turn it is upper bounded by $M^N$. In addition, $L \leq M$ and $2L \leq N^2 < 2(L+1) \leq 2(M+1) \leq 3M$,
\[
(LMN)^{12N} \leq \left(M^2 \sqrt{3N}\right)^{12N} =3^{6N} M^{\frac{5}{2} 12N} \leq M^{31N}.
\]
We can gather all powers in terms of $M$, so for all $z$ with $0 < |z| < r = \frac{1 + |q|}{2}$,
\begin{equation}
|F(z)| \leq |z|^M M^{31N}. \label{EqUpperbound}
\end{equation}

\section{Lower bound on prime powers $|F(q^P)|$}
Consider a prime number $P$ such that $F(q^P) \not= 0.$ Remark that such a $P$ exists, for that arguing over closed disks on $\D_{|q|} \subset\D$, $F \not\equiv 0$ is a holomorphic function, so it can only have finitely many zeros. Set now  $\al = A(q^P, J(q^P))$, then $F(q^P) = \Delta(q^P)^{2N}\al$. First observe that as $z^{-1}\Delta(z)$ is holomorphic and nonvanishing on $\D$, by the minimum principle we have an absolute constant $C_6 >0$ such that $|\Delta(q^P)| \geq C_6|q|^P$, hence 
\begin{equation}
|\Delta(q^P)|^{2N} \geq C_6^{2N} |q|^{P2N} \geq \exp(-(C_7\log(|q|^{-1})NP).\label{EqLowerBoundDelta}
\end{equation}
Let us focus now on $\al = A(q^P, J(q^P))$. As $q, J(q) \in \Qbar$, $q^P \in \Qbar$ and $J(q^P) \in \Qbar$ as $\Phi_P(J(q), J(q^P)) = 0$, therefore $\al \in \Qbar$ and we will use Liouville's inequality Lemma \ref{LemmaLiouvilleIneq}. We then need upper bounds on $\deg(\al)$ and $h(\al)$. As $\al \in \Q(q, J(q^P))$ and $J(q^P)$ belongs to an algebraic extension of $\Q(J(q))$ of degree bounded by $P+1,$ we have
\begin{align*}
\deg(\al) &\leq [\Q(q, J(q^P)): \Q] \leq [\Q(q):\Q] [\Q(J(q^P)):\Q] = \deg(q)\deg(J(q^P))\\
& \leq (P+1) \deg(q) \deg(J(q)).
\end{align*}
For the height, by  Lemma \ref{LemaHeightEvalPol},
\[
h(A(q^P, J(q^P)) \leq \log L(A) + Nh(q^P) + Nh(J(q^P)) = \log L(A) + NPh(q) + Nh(J(q^P)).
\]
For $h(J(q^P)),$ the original proof of \cite{Stephanois} compared the $h(J(q^P))$ with $h(J(q))$ via Lemma \ref{LemaMahlerMeasOnARoot}, with the modular polynomial $\Phi_P$ and $J(q),$
\begin{gather*}
\deg(J(q^P)) h(J(q^P)) \leq \deg(J(q)) \left(\log L(\Phi_P) + \deg_x(\Phi_P) h(J(q)) \right) \\
\leq \deg(J(q)) \left( C'(P+1)\log P + (P+1)h(J(q)) \right).
\end{gather*}
Alternatively, for a more geometric version, we use Proposition \ref{PropHeightIsogeJInv}:
\[
h(J(q^P)) \leq 2h(J(q)) + 6\log(1 + P) + C_2.
\]
Therefore, \small
\begin{align*}
\deg(\al) h(\al) & \leq \deg(\al)\left(\log L(A) + NP h(q) +Nh(J(q^P))\right), \\
&\leq (P+1)\deg(J(q))\deg(q) \left(\log L(A)+ NP h(q) + N\left(2h(J(q)) + 6\log(1 + P) + C_2\right)  \right),\\
&\leq_{\exists C_7} (P+1)\deg(J(q))\deg(q) \left(\log L(A)+ NP h(q) + NC_7\log(P) + 2Nh(J(q)) \right)
\end{align*}
\normalsize and as $\log L(A) \leq \log(C_3^N (N^2/2)^{12N}) \leq 25N \log N$, if we set \\$C_8(q) = \deg(q) \deg(J(q))\max\{1, 25, C_7, h(q), h(J(q)) \}$, then
\[
\deg(\al) h(\al) \leq C_8(q)N(P+1)\left( P + \log N + \log P + 1 \right) \leq_{\exists C_9(q)} C_9(q)NP\left( P + \log N \right).
\]
Then Liouville's inequality $\log|\al| \geq -\deg(\al) h(\al)$ and our previous bound \eqref{EqLowerBoundDelta} allow us to conclude that exists $C_{10}(q)$ such that
\begin{equation}
|F(q^P)| \geq \exp \left( -C_{10}(q)NP\left( P + \log N \right)\right). \label{EqLowerBoundF(qP)}
\end{equation}
\section{Definition of $P$ and upper bound on $P$}
Combining \eqref{EqUpperbound} and \eqref{EqLowerBoundF(qP)}:
\begin{align*}
\exp\left(-C_{10}(q)NP\left( P + \log N \right) \right) & \leq |q|^{PM} M^{31N} = \exp \left( -\log(|q|^{-1})PM + 31N \log M \right),\\
\log(|q|^{-1})PM & \leq C_{11}(q) NP\left( P + \log N \right) + 31N\log M, \\
M & \leq_{\exists C_{12}(q)} C_{12}(q) N\left(P + \log P + \log N + \frac{31\log M}{P}\right),
\end{align*}
and as $N^2 \leq 2(L+1) \leq 2(M+1) \leq 3M$, it follows that $\log N \leq \log M$ and $N \leq \sqrt{3M}$ and there exists $C_{13}(q)$ such that
\begin{equation}
M \leq C_{13}(q) \sqrt{M}\left(P + \log M \right). \label{EqContradictionOnM}
\end{equation}
Observe that at this point we are aiming for a contradiction of the sort $M =O(\sqrt{M} \log M)$, and  we would have that as long as we may control $P$  well enough in terms of $N,M$. This is also the step where the strategy differs from the original Mahler's method.

The only condition that we impose on $P$ is that $F(q^P) \not= 0$. We said before that such an $P$ exists, so we may define now $P$ such that $F(q^p) = 0$ for all primes $p < P$. 

Set $r = \frac{1 + |q|}{2}$. We can think of the following standard construction of holomorphic function on $\D_r$, 
\[
\frac{F(z)}{z^M} \prod_{p \text{ prime, } p <P } \frac{1}{z-q^p},
\]
and we are simply going to change the factors $ \frac{1}{z-q^p}$ to the following inverse Blaschke-type product
\[
H(z) : =\frac{F(z)}{z^M}  \prod_{p \text{ prime, } p < P} \frac{r^2 - z\bar{q}^p}{r(z-q^p)}.
\]
We consider these modified factors because at the boundary $\{ |z|= r \}$,\footnote{this argument is from \cite[paragraph under Equation (6)]{Netserenko}.}
\[
\lvert r^2 - \bar{q}^p z\rvert = \lvert r^2 - \bar{q}^p z\rvert \frac{|\bar{z}|}{r} = \frac{1}{r}\lvert r^2\bar{z} - \bar{q}^p r^2 \rvert = r\lvert z - \bar{q}^p \rvert = \lvert r(z - q^p) \rvert,
\]
so $\max_{|z|=r} |H(z)| =  r^{-M} |F(z)|$. By the maximum principle, for 
$z \in \D_{r}$:
\[
|H(z)| \leq r^{-M} |F(z)| \leq M^{31N},
\]
by \eqref{EqUpperbound}. On the other hand, $H(0) = d_0 \prod_{p \text{ prime, } p < P} \frac{-r}{q^p}$, with $d_0$ the first coefficients of the expansion of $G(z) = z^{-M}F(z)$. By construction, $d_0 \in \Z$ and it is not zero, so $|d_0| \geq 1$. Together with $|H(0)|\leq M^{31N},$ we have our bound involving $P, M, N$. We require to bound $\pi(P):=\sum_{p \text{ prime, } p < P} 1$ and $\sum_{p \text{ prime, } p < P} p$, and we use Proposition \ref{PropBoundsSumsPrimes}.
Therefore, we have:
\begin{align*}
r^{\pi(P)}|q|^{-\sum_{p \text{ prime, } p < P} p } &\leq M^{31N},\\
\left(\frac{1}{r}\right)^{-C_{14}P/\log P}\left(\frac{1}{|q|}\right)^{P^2/2\log P} & \leq M^{31N},\\
\frac{P^2}{2 \log P} \log (|q|^{-1}) - C_{14}\log(r^{-1})\frac{P}{\log P} & \leq 31N \log M,
\end{align*} 
and as $r \geq |q|$ it follow that $-\log r^{-1} \geq -\log |q|^{-1}$, so there exists $C_{15}>0$ such that $(P/2 - C_{14}) \geq (C_{15})^{-1}P$ and
\[
\frac{P^2}{2 \log P} \log (|q|^{-1}) - C_{14}\log(r^{-1})\frac{P}{\log P}  \geq \frac{P}{\log P} \log (|q|^{-1}) \left( \frac{P}{2} - C_{14} \right) \geq \frac{P^2}{\log P}\log (|q|^{-1}) \frac{1}{C_{15}},
\]
hence 
\begin{equation}
\frac{P^2}{\log P} \leq C_{15} \frac{1}{\log (|q|^{-1})} 31 N \log M. \label{EqBoundOnS}
\end{equation}

\section{Final Contradiction}
\label{sec:FinalContradiction}
In \eqref{EqBoundOnS}, we can take $\frac{P^2}{\log P} \geq P\sqrt{P}$, so the following arguments are easier. Using $N \leq \sqrt{3M}$, it follows
\[
P^{3/2} \leq \frac{C_{15}}{\log (|q|^{-1})}31N \log M \leq_{\exists C_{16}(q)} C_{16}(q) \sqrt{M}\log M.
\]
Therefore, together with \eqref{EqContradictionOnM}, there exists $C_{17}(q)$ such that
\begin{align*}
M & \leq C_{17}(q) \sqrt{M} \left( \log M + \left( \sqrt{M} \log M \right)^{2/3} \right) \leq_{\exists C_{18}(q)} C_{18}(q) \sqrt{M} \left( \sqrt{M}\log M \right)^{2/3},\\
M & \leq C_{18}(q) M^{5/6} (\log M)^{2/3},
\end{align*}
which is a contradiction for $M$ large enough. As $M \geq [N^2/2]$, and $N$ can be taken arbitrarily large,\footnote{recall that the only dependence of $N$ on $q$ is in Equation \eqref{EqDependNonq}.} we have our final contradiction.

\appendix
\section{Other bounds on $P$ that do not yield a contradiction}
This section is no part of the proof, but it is included with sights to generalizations for genus 2. 

The following bounds are weaker, of $P = O(N)$ and $P = O(N \log M)$, and one can see that in the previous section that they do not yield a contradiction. Actually, any bound of the form $P^{1 + \eps} = O(N \log M)$ for any $\eps >0$ gives the contradiction, but not for $\eps = 0$.

\subsection{Jensen's formula}
We have a holomorphic function $G(z) = z^{-M}F(z)$ such that on $D_r$ it follows $|G(z)| \leq M^{31N}$, and we want to derive information of the zeroes form the growth of $G$. Jensen's formula (\cite[Lemma 6.1]{MontVauMNT} with $r' = |q| < r = \frac{1 + |q|}{2}$) implies that the number of zeros of $G$ in $\{|z| < |q| \}$ does not exceed
\[
\frac{1}{\log (r/r')} \frac{1}{\underbrace{|G(0)|}_{\geq 1}} \log(M^{31N})= O_q(N\log M),
\]
as $G(0) = d_0 \in \Z$. In particular, $P = O_q(N\log M)$. From this perspective, one can understand the strategy in the last section as an strengthening of the bound given by Jensen's formula, by using that the zeros one is counting are in a geometric sequence $\{q^p\}$ for $p$ prime.

\subsection{A purely algebraic argument}
\label{ssec:PurelyAlgebraic}
This argument shows well-definedness of $P$ without using complex analysis. Suppose $A(q^p, J(q^p) = 0$ for $p$ prime, with $(q,J(q)) \in \Qbar^2$. This realizes $J(q^p)$ as a root of a polynomial in $\Q(q^p)[X] \subset \Q(q)[X]$ of degree $N$, so
\[
\deg(J(q^p)) = [\Q(J(q^p)): \Q] \leq [\Q(J(q^p): \Q(q)] [\Q(q): \Q] \leq \deg(q) N,
\]
but this bound is independent of $p$, which is absurd: by Proposition \ref{PropBertrand}, for $p$ large enough in terms of $J(q)$, 
\[
\frac{p-1}{3} \leq [\Q(J(q),J(q^p)) : \Q(J(q))].
\]
This lower bound is only necessary for CM elliptic curves, for the non-CM for $p$ large enough,\footnote{larger than the prime in Serre's uniformity theorem.} the degree is generically $p+1$.

Hence, as $\deg(J(q^p)) = [\Q(J(q^p)): \Q] \geq [\Q(J(q))(J(q^p)): \Q(J(q))] = [\Q(J(q),J(q^p)) : \Q(J(q))]$, it follows that
\[
\frac{p-1}{3} \leq \deg(q)N,
\]
so taking $P$ as the smallest $p$ such that this bound does not hold (and among the large enough primes for Proposition \ref{PropBertrand} to apply, but this does not depend on $N$), we have $A(q^P, J(q^P)) \not= 0$ and a bound $P = O_{q}(N)$. 
\printbibliography
\end{document}